\documentclass[12pt]{amsart}
\usepackage[all]{xypic}
\usepackage{amsmath,graphics}
\usepackage{amsfonts,amssymb}

\theoremstyle{plain}
\newtheorem*{theorem*}{Theorem}
\newtheorem*{lemma*} {Lemma}
\newtheorem*{corollary*} {Corollary}
\newtheorem*{proposition*} {Proposition}
\newtheorem{theorem}{Theorem}[section]
\newtheorem{lemma}[theorem]{Lemma}
\newtheorem{corollary}[theorem]{Corollary}

\theoremstyle{remark}
\newtheorem*{remark}{Remark}
\newtheorem*{remarks}{Remarks}

\theoremstyle{definition}
\newtheorem{definition}[theorem]{Definition}

\def\eps{\epsilon}

\def\old{\ol{\delta}}

\def\gl{\mbox{GL}}
\def\Q{\Bbb{Q}}
\def\k{\Bbb{K}}
\def\K{\Bbb{K}}

\def\id{\mbox{id}}

\def\Z{\Bbb{Z}}
\def\R{\Bbb{R}}

\def\N{\Bbb{N}}
\def\l{\lambda}
\def\part{\partial}

\def\a{\alpha}

\def\g{\gamma}

\def\bp{\begin{pmatrix}}

\def\sm{\setminus}
\def\ep{\end{pmatrix}}
\def\bn{\begin{enumerate}}

\def\rank{\mbox{rank}}

\def\en{\end{enumerate}}
\def\ba{\begin{array}}
\def\ea{\end{array}}

\def\S{\Sigma}

\def\v{\alpha}

\def\ti{\tilde}
\def\lk{\mbox{lk}}

\def\fr12{\frac{1}{2}}

\def\ker{\mbox{Ker}}

\def\hom{\mbox{Hom}}

\def\kt{\K\tpm}

\def\deg{\mbox{deg}}

\def\gcd{\mbox{gcd}}
\def\v{\varphi}
\def\tpm{[t^{\pm 1}]}

\def\be{\begin{equation}}
\def\ee{\end{equation}}

\def\F{\Bbb{F}}
\def\k{\Bbb{K}}
\def\K{\Bbb{K}}

\def\kat{\k[t^{\pm 1}]}
\def\ktfield{\k(t)}

\def\kat{\k[t^{\pm 1}]}

\def\ol{\overline}

\def\K{\Bbb{K}}

\def\tnphi{||\phi||_T}

\def\modd{\, mod \, }

\def\cmtbf#1{} \def\cmt#1{}

\begin{document}

\title{The parity of the Cochran--Harvey invariants of 3--manifolds} \author{Stefan Friedl
and Taehee Kim}
\date{\today} \address{Rice University, Houston, Texas, 77005-1892}
\email{friedl@math.rice.edu}\address{ Department of Mathematics, Konkuk University,
Hwayang-dong, Gwangjin-gu, Seoul 143-701, Korea} \email{tkim@rice.edu}
\def\subjclassname{\textup{2000} Mathematics Subject Classification}
\expandafter\let\csname subjclassname@1991\endcsname=\subjclassname \expandafter\let\csname
subjclassname@2000\endcsname=\subjclassname 
\date{\today}
\begin{abstract}
Given a finitely presented group $G$ and an epimorphism $\phi:G\to \Z$  Cochran and Harvey defined
a sequence of invariants $\old_n(G,\phi)\in \N_0, n\in \N_0$, which can be viewed as the degrees of
higher--order Alexander polynomials. Cochran and Harvey showed that (up to a minor modification)
this is a never decreasing sequence of numbers if $G$ is the fundamental group of a 3--manifold
with empty or toroidal boundary. Furthermore they showed that these invariants give lower bounds on
the Thurston norm.

Using a certain Cohn localization and the duality of Reidemeister torsion we show that for a
fundamental group of a 3--manifold any jump in the sequence is necessarily even. This answers in
particular a question of Cochran. Furthermore using results of Turaev we show that under a mild
extra hypothesis the parity of the Cochran--Harvey invariant agrees with the parity of the Thurston
norm.
 \end{abstract}
\maketitle

%Hi Taehee,
%
%I mostly followed the referee's very reasonable suggestions. In particular I tried to shorten some
%of our long sentences. The only other modification is that I was trying to make sure our theorems
%are correct even in silly cases like $\Delta_K=1$, $M=S^1\times D^2$ etc. Perhaps when you read the
%papers you can keep an eye on such extreme cases.

%Ideas, questions:
%\bn
%\item I presume that the parity results are false in general for 2--complexes with $\chi(X)=0$.
%\item it shouldn't be so difficult to define monicness of element in
%$\ktfield^*/[\ktfield^*,\ktfield^*]$. Of course not quite clear how this can detect satellite
%knots.
%\en

%=====================================
\section{Introduction} \label{section:intro}
Let $G$ be a finitely presented group and $\phi:G\to \Z$ an
epimorphism. Following \cite{Ha05,Ha06,Co04} we define for $n\geq
0$
\[ \delta_n(G,\phi):=\deg(\Delta_n(G,\phi)(t)),\]
where $\Delta_n(G,\phi)(t)$ denotes the $n$--th higher--order Alexander polynomial of $G$
corresponding to $\phi$. These polynomials are elements in $\K_n\tpm$, where $\K_n\tpm$ is a skew
Laurent polynomial ring which is associated to the group ring $\Z[G/G_r^{(n+1)}]$ and to the
homomorphism $\phi$. Here $G_r^{(n+1)}$ denotes the $(n+1)$--st term of Harvey's rational derived
series. We adopt the convention that $\deg(0)=-\infty$. We give the precise definitions in Section
\ref{section:admissible}. Note that in \cite{Ha06}  these invariants are denoted by
$\ol{\delta}_n(G,\phi)$.

These invariants show an interesting behavior for the fundamental groups of 3--manifolds. (Here by
a 3--manifold we always mean an oriented, connected and compact 3--manifold.) Let $G$ be the
fundamental group of a 3--manifold $M$ with empty or toroidal boundary. Then Harvey \cite{Ha06}
showed (cf. also \cite{Co04} and \cite{Fr06}) that if  $\phi:G\to \Z$ is an epimorphism such that
$\delta_0(G,\phi)\ne -\infty$, then
\[ \ba{rccl} \delta_0(G,\phi)-1-b_3(M) &\leq &\delta_1(G,\phi)  \leq \delta_2(G,\phi) \dots
&\mbox{ if $b_1(M)=1$, and}\\
\delta_0(G,\phi)& \leq &\delta_1(G,\phi)  \leq \delta_2(G,\phi) \dots &\mbox{ if $b_1(M)>1$.}\ea
\] Here $b_i(M)$ denotes the $i$--th Betti number of $M$. Recall that $b_3(M)=1$ if $M$ is closed and
$b_3(M)=0$ if $M$ has boundary. This result can be viewed as a subtle obstruction to a group being
a 3--manifold group. Note that a similar result holds for groups of deficiency one (cf. \cite{Ha06}
or \cite{Fr06}).

If $\delta_0(G,\phi)=0$ and $b_1(M)=1$, then the proof of Lemma \ref{lem:cyclic} shows that
$\delta_n(G,\phi)=0$ for all $n>0$. We can and will therefore restrict ourselves to the cases
$\delta_0(G,\phi)>0$ or $b_1(M)\geq 2$.
%\footnote{The next theorem would be wrong otherwise, e.g.
%if $G=\pi_1(S^3\sm K), \Delta_K(t)=1$, since then $\delta_i(G)=\delta_0(G)$}

The following theorem strengthens Harvey's result:
\begin{theorem}\label{mainthm1}
Let $G$ be the fundamental group of a 3--manifold with empty or toroidal boundary. Let $\phi:G\to
\Z$ be an epimorphism such that $\delta_0(G,\phi)\ne -\infty$. Then any jump in the sequence of
inequalities
\[ \ba{rccl} \delta_0(G,\phi)-1-b_3(M) &\leq &\delta_1(G,\phi)  \leq \delta_2(G,\phi) \dots \,
&\mbox{(if $b_1(M)=1$ and $\delta_0(G,\phi)>0$)}\\
\delta_0(G,\phi)& \leq &\delta_1(G,\phi)  \leq \delta_2(G,\phi) \dots \, &\mbox{(if $b_1(M)>1$)}\ea
\]
has to be even.
\end{theorem}
Note that this theorem is not known to hold for groups of deficiency one in general. This theorem
in particular answers a question of Cochran (cf. \cite[Question (5) in Section 14]{Co04}).
Furthermore, together with realization results of Cochran, this completely classifies the sequences
of numbers which appear as $\delta_n$--invariants for knots (see Theorem \ref{thm:paritydelta}).

The key difficulty in relating $\delta_{k+1}(G,\phi)$ and $\delta_{k}(G,\phi)$ is that there are no
non--trivial homomorphisms between the skew Laurent polynomial rings $\K_{k+1}\tpm$ and $\K_k\tpm$.
We circumvent this problem by studying an appropriate Cohn localization which maps non--trivially
to $\K_{k+1}\tpm$ and to $\K_k\tpm$, and by using the functoriality  of Reidemeister torsion.
Another major ingredient is a duality theorem for Reidemeister torsion.
\\

The Cochran--Harvey invariants have a close relationship to the geometry of 3--manifolds. Let $M$
be 3--manifold $M$ and $\phi\in H^1(M;\Z)$. Then  the \emph{Thurston norm} of $\phi$ is defined as
 \[
||\phi||_{T}=\min \{ -\chi(\hat{S})\, | \, S \subset M \mbox{ properly embedded surface dual to
}\phi\}
\] where $\hat{S}$ denotes the result of discarding all connected components of $S$ with positive Euler
characteristic. This extends to a seminorm on $H^1(M;\R)$. We refer to \cite{Th86} for details. We
henceforth identify $H^1(M;\Z)$ with $\hom(H_1(M;\Z),\Z)$ and $\hom(\pi_1(M), \Z)$. Let $M$ be
3--manifold with empty or toroidal boundary and let $G=\pi_1(M)$. In \cite{Ha05,Ha06} Harvey (cf.
also \cite{Co04}, \cite{Tu02} and \cite{Fr06}) showed that for an epimorphism  $\phi:G\to \Z$ with
$\delta_0(G,\phi)\ne -\infty$ we have
\begin{equation}\label{equ:lowerbound} \ba{rccl}
\delta_0(G,\phi)-1-b_3(M) &\leq &\delta_1(G,\phi)  \leq \delta_2(G,\phi) \dots \leq ||\phi||_T
&\mbox{(if $b_1(M)=1)$}\\
\delta_0(G,\phi)& \leq &\delta_1(G,\phi)  \leq \delta_2(G,\phi) \dots \leq ||\phi||_T &\mbox{(if
$b_1(M)>1$).}\ea
\end{equation}
Furthermore if $(M,\phi)$ fibers over $S^1$ and if $M\ne S^1\times S^2, S^1\times D^2$, then all
inequalities are in fact equalities. Using classical results of Turaev \cite{Tu86} we can
strengthen this geometric interpretation as follows.

\begin{theorem}\label{mainthm2}
Let $M$ be a  3--manifold which is either closed or the exterior of a link in $S^3$. Let
$G:=\pi_1(M)$ and let $\phi:G\to \Z$ be an epimorphism such that $\delta_0(G,\phi)\ne -\infty$,
then
\[ \ba{rccl} \delta_0(G,\phi)-1-b_3(M) &\hspace{-0.2cm}\equiv &\hspace{-0.2cm}\delta_1(G,\phi)
\equiv \delta_2(G,\phi) \dots \equiv ||\phi||_T \modd 2
&\hspace{-0.2cm}\mbox{(if $b_1(M)=1$, $\delta_0(G,\phi)\geq 1$)}\\
\delta_0(G,\phi)&\hspace{-0.2cm} \equiv &\hspace{-0.2cm}\delta_1(G,\phi) \equiv \delta_2(G,\phi)
\dots \equiv ||\phi||_T  \modd 2  &\hspace{-0.2cm}\mbox{(if $b_1(M)>1$).}\ea
\]
%\[ \delta_0(G,\phi)-1-b_3(M) \equiv \delta_1(G,\phi) \equiv \delta_2(G,\phi) \dots \equiv ||\phi||_T\, \modd 2.\]
\end{theorem}

In a forthcoming paper the first author and Shelly Harvey \cite{FH06} will show that  for a
3--manifold group $G$ the map $\phi\mapsto \delta_n(G,\phi)$ defines  a seminorm on $H^1(G;\R)$.
This further strengthens the close relationship between the Cochran--Harvey invariants and the
Thurston norm.

The paper is organized as follows. In Section \ref{section:rt} we recall several facts about
Reidemeister torsion and fix some notation. We also define the Alexander polynomials and recall
the relationship to Reidemeister torsion. In Section \ref{section:admissible} we recall several
definitions from \cite{Ha06} which will then allow us to give in Section \ref{section:results}
precise formulations of (slight generalizations of) Theorems \ref{mainthm1} and \ref{mainthm2}.
We conclude with the proofs of our results in Section \ref{section:proof}.
\\

{ \bf Acknowledgment:} We would like to thank the anonymous referee for helpful
 suggestions for improving the clarity of the paper and for pointing out a
 mistake in an earlier version.

%=====================================
\section{Reidemeister torsion and Alexander polynomials} \label{section:rt}

%=====================================
\subsection{Reidemeister torsion of manifolds} Throughout this paper we will only
consider associative rings $R$ with $1\ne 0$ with the property that if $r\ne s$, then $R^r$ is
not isomorphic to $R^s$. For such a ring $R$, let $\gl(R):=\underset{\rightarrow}{\lim} \,
\gl(R,n)$ and $K_1(R):=\gl(R)/[\gl(R),\gl(R)]$.
%In the direct system we have the maps $\gl(R,n)\to
%\gl(R,n+1)$ given by $A \mapsto \bp A &0\\0&1\ep $.
%Denote by $E(R)$ the subgroup given by all elementary matrices, i.e. all matrices
%$E=(e_{ij})$ of the form $e_{ii}=1$ for all $i$ and $e_{ij}=0$ except for one pair $(i,j)$. Then
%$E(R)=[\gl(R),\gl(R)]$
%We define .
In particular $K_1(R)$ is an abelian group. If $\K$ is a (skew) field, then the Dieudonn\'e
determinant gives an isomorphism $\det:K_1(\K)\to \K^\times/[\K^\times,\K^\times]$ where
$\K^\times:=\K\sm \{0\}$. For more details we refer to \cite{Ro94} or \cite{Tu01}.

Let $X$ be a CW--complex. In this paper by a CW--complex we will always mean a  connected
CW--complex. Denote the universal cover of $X$ by $\ti{X}$. We view the chain complex $C_*(\ti{X})$
as a right $\Z[\pi_1(X)]$--module via deck transformations. Let $\v:\Z[\pi_1(X)]\to R$ be a
homomorphism. This equips $R$ with a left $\Z[\pi_1(X)]$--module structure. We can therefore
consider the right $R$--module chain complex $C_*(X;R):=C_*(\ti{X})\otimes_{\Z[\pi_1(X)]} R$. We
denote the homology of $C_*(X;R)$ by $H_i(X;R)$.

Now assume that $X$ is in fact a finite CW--complex. If $H_*(X;R) = 0$, then the
\emph{Reidemeister torsion} $\tau(X,R)\in K_1(R)/\pm \v(\pi_1(X))$ is defined (cf. \cite{Tu01}).
We will write $K_1(R)/\pm \pi_1(X)$ for $K_1(R)/\pm \v(\pi_1(X))$ for convenience, if $\v$ is
clear from the context. If $H_*(X;R) \ne 0$, then we write $\tau(X,R) :=0$.
%The
%Reidemeister torsion of $(X,R)$ depends on the choice of the homomorphism $\v$, but we suppress
%$\v$ in the notation for convenience. The homomorphism $\v$ is normally clear from the context.
It is known that $\tau(X,R)$ only depends on the homeomorphism type of $X$. Hence we can define
$\tau(M,R)$ for a 3--manifold $M$ by picking any finite CW--structure for $M$. We refer to the
excellent book of Turaev \cite{Tu01} for filling in the details.

%=====================================
\subsection{Alexander polynomials}\label{sec:phicompatible}
%The following algebraic setup allows us to define non--commutative
%Alexander polynomials. First

Let $\k$ be a (skew) field and $\g:\k\to \k$ a ring homomorphism. Denote by $\kat$ the skew Laurent
polynomial ring over $\k$ (associated to $\g$). The elements in $\kat$ are formal sums
$\sum_{i=-r}^s a_it^i$ with $a_i\in \k$ and multiplication in $\kat$ is induced from the rule
$t^ia=\g^i(a)t^i$ for any $a\in \k$. It is known that $\kat$ embeds in its (skew) quotient field
(denoted by $\ktfield$)  which is flat over $\kat$ (cf. e.g. \cite{DLMSY03}).

Let $X$ be a CW--complex with finitely many cells in dimension 1 and 2.
% Using $\v$ we equip $\kat$ with a (left) $\Z[\pi_1(X)]$--module structure.
%Theorem \ref{thm:tfa} can easily be extended to show that $\kat$ also has an Ore localization
%(cf. also \cite{DLMSY03}) which is flat over $\kt$.
Given a ring homomorphism $\v:\Z[\pi_1(X)]\to \kt$ to a skew Laurent polynomial ring we  consider
the finitely presented $\kat$--module $H_1(X;\kat)$.
% := H_i(C_*(\ti{X})\otimes_{\Z[\pi_1(X)]} \kat)$ where $\ti{X}$ is the
%universal cover of $X$. We drop $\v$ in the notation if it is clear from the context.
This module was studied in \cite{Co04,Ha05,Tu02,LM05}. Since $\kat$ is a principal ideal domain
(PID) (cf. \cite[Proposition~4.5]{Co04}) we can decompose
\[ H_1(X;\kat)\cong \bigoplus_{i=1}^l
\kat/(p_i(t))\] for $p_i(t)\in \kat$, $1 \le i \le l$. We define
$\Delta^\v(t):=\prod_{i=1}^lp_i(t) \in \kat$. $\Delta^\v(t)$ is called the  {\em Alexander
polynomial} of $(X,\v)$ and it is well-defined up to a certain indeterminacy which is determined
in \cite[Theorem~3.1]{Fr06} (cf. also \cite[p.~367]{Co04} and \cite{Co85}). If $X=K(G,1)$ for a
group $G$, then we also call $\Delta^\v(t)$ the Alexander polynomial of $(G,\v)$.
%We refer to  for a discussion of the indeterminacy of $\Delta^\v_{\phi}(t)$.
%We drop the subscript $\phi$ when $\phi$ is clear from the context.
% In fact
%the only known invariant which can be extracted from the non--commutative Alexander polynomials
%are their degrees.

%We call $\phi \in H^1(X;\Z)$ {\em primitive} if the corresponding map $\phi:H_1(X;\Z)\to \Z$ is
%surjective. In \cite{Fr06} we prove the following lemma.
%
%\begin{lemma} \label{lem:alex03}
%Let $M$ be a 3--manifold, $\phi \in H^1(M;\Z)$ primitive. Let  $\v:\pi_1(X)\to \kat$ be a
%$\phi$--compatible homomorphism.
% \bn
% \item If $\im(\v(\pi_1(X)))\subset \kat$ is cyclic, then
%$\Delta^\v_{0}(t)=at-1$ for some $a\in \k$. Otherwise $\Delta^\v_{0}(t)=1$.
% \item Assume that $\Delta^\v_{1}(t)\ne 0$. If $M$ has boundary,
%then $\Delta^\v_{2}(t)=1$, otherwise $\Delta^\v_{2}(t)=\Delta^\v_{0}(t)$.
%\en
%\end{lemma}
%
%Based on ideas of Turaev we prove the following theorem in \cite{Fr06}, relating the
%non--commutative Alexander polynomials to Reidemeister torsion.

%The following theorem gives the relationship between $\tau(X,\ktfield)$ and Alexander
%polynomials.
%\begin{theorem}\label{thm:torsionalex}\cite[Theorem 1.1]{Fr06}
%Let $M$ be a 3--manifold whose boundary is empty or consists of
%tori. Let $\phi\in H^1(M;\Z)$ be non--trivial and let
%$\v:\Z[\pi_1(M)] \to \kt$ be a $\phi$--compatible homomorphism.
%Then $\tau(M,\ktfield) \ne 0$ if and only if $\Delta^\v_1(t)\ne
%0$. Furthermore if $\Delta^\v_1(t)\ne 0$, then
%\[ \tau(M,\ktfield)=\Delta^\v_0(t)^{-1}\Delta^\v_1(t)\Delta^\v_2(t)^{-1} \in K_1(\ktfield)/K_1(\kt).\]
%\end{theorem}

Let $f(t)\in \kat$. If $f(t)=0$ then we write $\deg(f(t))=-\infty$. Otherwise, for
$f(t)=\sum_{i=m}^n a_it^i\in \kat$ with $a_m\ne 0, a_n \ne 0$ we define $\deg(f(t)):=n-m$. Note
that $\deg(\Delta^\v(t))$ is well--defined (cf. \cite{Co04}). In fact we have the following
interpretation of the degree of the Alexander polynomial:
\be \label{equ:degdim} \deg(\Delta^\v(t))=\dim_{\K}(H_1(X;\kat)).\ee

Given $\v:\Z[\pi_1(X)] \to \kat$ we also consider the induced map $\Z[\pi_1(X)] \to \kat\to
\ktfield$. In particular we can consider $\tau(X,\ktfield)\in K_1(\ktfield)/\pm \pi_1(X) \, \cup
\{0\}$. Note that setting $\deg(f(t)g(t)^{-1})=\deg(f(t))-\deg(g(t))$ for $f(t),g(t)\in \kt$ we can
extend $\deg:\kt\sm \{0\} \to \Z$ to a homomorphism $\deg:\ktfield^\times\to \Z$ which in turn
induces (using the Dieudonn\'e determinant) a homomorphism $\deg:K_1(\ktfield)\to \Z$. Note that a
ring homomorphism $ \Z[\pi_1(X)] \xrightarrow{\v} \kat$ sends $g\in \pi_1(X)$ to an invertible
element in $\kat$, i.e. to an element of the form $kt^l$, whose degree is zero. We therefore see
that $\deg:K_1(\ktfield)\to \Z$  passes to $\deg : K_1(\ktfield)/\pm\pi_1(X)\, \cup \{0\} \to
\Z\cup \{-\infty\}$.

Let $\phi \in H^1(X;\Z)$. We say $\phi\in H^1(X;\Z)$ is  \emph{primitive} if the induced
homomorphism $\phi : \pi_1(X) \to \Z$ is surjective.
 Following
Turaev \cite{Tu02} we call a homomorphism $\v:\Z[\pi_1(X)] \to \kat$ \emph{$\phi$--compatible}
if for any $g\in \pi_1(X)$ we have $\v(g)=kt^{\phi(g)}$ for some $k\in \k$.

%Combining Lemmas 3.1 and 3.2 in \cite{Fr06} with Theorem
%\ref{thm:torsionalex} we get the following corollary.

\begin{theorem} \label{thm:taudelta}
Let $M$ be a 3--manifold with empty or toroidal boundary. Let $\phi \in H^1(M;\Z)$ primitive and
let $\v:\Z[\pi_1(M)]\to \kat$ be a $\phi$--compatible homomorphism.
%Suppose the same hypotheses as in Theorem \ref{thm:torsionalex}.
If $\Delta^\v(t)\ne 0$, then $\tau(M,\ktfield)\ne 0$. Furthermore, if $\v(\pi_1(M))\subset \kat$
is cyclic, then
\[  \deg(\tau(M,\ktfield))=\deg(\Delta^\v(t))-(1+b_3(M)), \]
otherwise
 \[  \deg(\tau(M,\ktfield))= \deg(\Delta^\v(t)). \]
\end{theorem}

\begin{proof}
Assume that $\Delta^\v(t)\ne 0$, then $\tau(M,\ktfield)\ne 0$ by \cite[Theorem~1.2]{Fr06}. It
follows from \cite[Theorem 1.1]{Fr06} that $\deg(\det(\tau(M,\ktfield)))$ is the alternating sum of
the degrees of the Alexander polynomials corresponding to $H_i(M;\kt), i=0,1,2$. A direct
computation of $H_0(M;\kt)$ and $H_2(M;\kt)$ as in \cite[Lemmas 4.3 and 4.4]{Fr06} concludes the
proof.
\end{proof}

%=============================================================
\section{Admissible pairs and triples} \label{section:admissible}

%==================================
%\subsection{Skew fields of group rings} \label{section:skewfield}
%A group $G$ is called locally indicable if for every finitely
%generated subgroup $U\subset G$ there exists a non--trivial
%homomorphism $U\to \Z$.

It is known that for a locally indicable and amenable group $G$ the group ring $\Z[G]$ embeds in
its (skew) quotient field which is flat over $\Z[G]$ (cf. e.g. \cite{Hi40},
\cite[Corollary~6.3]{DLMSY03} and  \cite[p.~99]{Ra98}).
 We denote the quotient field by $\K(G)$. For instance, the
poly-torsion-free-abelian (henceforth PTFA) groups studied in \cite{COT03,Co04,Ha05} are locally
indicable and amenable (cf. e.g. \cite{St74}).

%\begin{theorem}\label{thm:tfa}
%Let $G$ be a locally indicable and amenable group and let $R$ be a subring of $\C$.
%\bn
%\item  $R[G]$  is an Ore domain, in particular it embeds in its classical right ring of
%quotients $\K(G)$.
%\item $\K(G)$ is flat over $R[G]$.
%\en
%\end{theorem}

%It follows from \cite{Hi40}  that $R[G]$ has no zero divisors. The first part now follows from
%\cite{Ta57} or \cite[Corollary~6.3]{DLMSY03}. The second part is a well--known property of Ore
%localizations (cf. e.g. \cite[p.~99]{Ra98}).
%
%A group $G$ is called poly--torsion--free--abelian (PTFA) if there exists a filtration
%\[ 1=G_0 \subset G_1\subset \dots \subset G_{n-1}\subset G_n=G \]
%such that $G_{i}/G_{i-1}$ is torsion free abelian. It is well--known that PTFA groups are
%amenable and locally indicable (cf. \cite{St74}). The group rings of PTFA groups played an
%important role in \cite{COT03}, \cite{Co04} and \cite{Ha05}.

%We will always consider the group ring $\Z[G]$ together with the
%involution given by $\ol{g}:=g^{-1}$ for $g\in G$. Note that this
%involution extends to involutions on $\K(G)$ and $K_1(\K(G))$.

%==================================
%\subsection{Admissible pairs} \label{section:admissible}

\begin{definition} Let $\pi$ be a group, $\phi:\pi \to \Z$ an epimorphism and
$\varphi:\pi\to G$ an epimorphism to a locally indicable and amenable group $G$ such that there
exists a map $\phi_G:G\to \Z$ (which is necessarily unique) such that
 \[   \xymatrix {
 \pi \ar[dr]_{\phi} \ar[r]^{\varphi} &G \ar[d]^{\phi_G}\\& \Z }
     \]
 commutes. Following \cite[Definition~1.4]{Ha06} we call $(\varphi,\phi)$
an {\em admissible pair for $\pi$}.
\end{definition}

Let $(\varphi:\pi\to G,\phi)$ be an admissible pair for $\pi$. Let $G' := \ker\{\phi_G:G\to \Z\}$.
Then $G'$ is still a locally indicable and amenable group, hence embeds in its (skew) quotient
field $\K(G')$. We pick an element $\mu \in G$ such that $\phi_G(\mu)=1$ and define $\g:\K(G')\to
\K(G')$ to be the homomorphism induced by $\g(g):=\mu g\mu^{-1}$ for $g \in G'$. Then we obtain a
Laurent polynomial ring $\K(G')\tpm$ with multiplication induced via $\g$ as described in the first
paragraph in Section \ref{sec:phicompatible}. We point out that different choices of $\mu$ give
isomorphic rings. Consider the ring homomorphism
\[ \ba{rcl} \Z[G]&\to& \K(G')\tpm\\[1mm]
    \sum n_gg &\mapsto& \sum n_gg\mu^{-\phi_G(g)}t^{\phi_G(g)}.\ea
\]
Clearly the induced homomorphism $\v:\Z[\pi]\to \Z[G]\to \K(G')\tpm$ is $\phi$--compatible.  Also
the map $f(t)g(t)^{-1}\mapsto f(\mu)g(\mu)^{-1}$ defines isomorphisms $\K(G')(t) \to \K(G)$ and
$\Z[G']\tpm \to \Z[G]$ (cf. also \cite[p.~364]{Co04}).

%Note that in particular $(\varphi_i,\phi), i=1,2$ are admissible
%pairs for $\pi$.

%We write $\tau(M,\K(G')(t))\in K_1(\K(G')(t))\cup \{0\}$ for the Reidemeister torsion. Given
%only a homomorphism $\phi:\pi\to \Z$ we can consider $\tau(M,\Q(t))\in K_1(\Q(t))\cup \{0\}$.

An important family of examples of admissible pairs is provided by Harvey's rational derived
series of a group $\pi$ (cf. \cite[Section~3]{Ha05}) which is defined as follows: Let
$\pi_r^{(0)}:=\pi$ and define inductively for $n\geq 1$
\[ \pi_r^{(n)}:=\big\{ g\in \pi_r^{(n-1)} | \, g^d \in \big[\pi_r^{(n-1)},\pi_r^{(n-1)}\big] \mbox{ for some }d\in \Z \sm \{0\} \big\}.\]
Note that $\pi_r^{(n-1)}/\pi_r^{(n)}$ is isomorphic to
$\big(\pi_r^{(n-1)}/\big[\pi_r^{(n-1)},\pi_r^{(n-1)}\big]\big)/\mbox{$\Z$--torsion}$. Hence it
is torsion--free abelian.
 By \cite[Corollary~3.6]{Ha05} the quotients $\pi/\pi_r^{(n)}$ are PTFA groups
for any $\pi$ and any $n$. If $\phi:\pi\to \Z$ is an epimorphism, then $(\pi\to
\pi/\pi_r^{(n)},\phi)$ is an admissible pair for $\pi$ for any $n>0$. Let $\pi' := \ker\{\phi : \pi
\to \Z\}$. Note that $\ker\{\pi/\pi_r^{(n+1)} \to \Z\} = \pi'/\pi_r^{(n+1)}$. We write
$\K_n:=\K\big(\pi'/\pi_r^{(n+1)}\big)$.
%In particular $\K_0=\Q$.
%Clearly for any epimorphism $\phi : \pi \to \Z$ and  any $n \geq 1$ we  get an admissible
%pair$(\pi\to \pi/\pi_r^{(n)},\phi)$ for $\pi$.

Given a finitely presented group $\pi$ and a homomorphism $\phi:\pi\to \Z$ we now define
$\Delta_n(\pi,\phi)(t)\in \k_n\tpm$ to be the Alexander polynomial of $(K(\pi,1),\Z[\pi]\to
\K(\pi'/\pi_r^{(n+1)})\tpm)$. This is called the \emph{$n$--th higher--order Alexander
polynomial of $\pi$ corresponding to $\phi$}. Note that this definition only works since we can
assume that $K(\pi,1)$ has finitely many cells in dimensions 1 and 2 since $\pi$ is finitely
presented.

\begin{remark}
Harvey \cite{Ha06} and Cochran \cite[Definition~5.3]{Co04} show the following equality: \be
\label{def:deltantwo}
\delta_n(\pi,\phi)=\dim_{\K(\pi'/\pi_r^{(n+1)})}(\pi_r^{(n+1)}/[\pi_r^{(n+1)},\pi_r^{(n+1)}]
\otimes_{\Z[\pi'/\pi_r^{(n+1)}]}\K(\pi'/\pi_r^{(n+1)})).\ee Here $\pi'/\pi_r^{(n+2)}$ acts on
the abelian group $\pi^{(n+1)}/\pi_r^{(n+2)}$ via conjugation. This equality can be easily
deduced from Lemma \ref{lemma:k1}, Equality (\ref{equ:degdim}) and the fact that Ore
localizations are flat.
\end{remark}

The next lemma follows immediately from the observation that  we can view $K(\pi,1)$ as the result
of adding $k$--cells with $k>2$ to $X$. It provides the link between group theory and topology.

\begin{lemma}\label{lemma:k1}
%\footnote{This lemma is new, it's important to fix the precise
%relationship between the theorems in Section 4 and the introduction. In order to prove it I
%added the above remark}
Let $X$ be a finite CW--complex with $\pi := \pi_1(X)$ and $\phi:\pi\to
\Z$ a non--trivial homomorphism. Then the Alexander polynomials of $(X,\Z[\pi]\to
\K(\pi'/\pi_r^{(n+1)})\tpm)$ and $(K(\pi,1),\pi\to \K(\pi'/\pi_r^{(n+1)})\tpm)$ agree.
\end{lemma}

In our applications of Theorem \ref{thm:taudelta} we will need the following basic lemma.
\begin{lemma} \label{lem:cyclic}
Let $\pi$ be a group and $n\geq 1$. Then $\pi/\pi_r^{(n)}$ is cyclic if and only if
\bn
\item $n=1$ and $\rank(H_1(\pi))=1$, or
\item $n\geq 2$, $\rank(H_1(\pi))=1$ and $\delta_0(\pi,\phi)=0$ where $\phi:\pi\to \Z$ is either of
the two epimorphisms.
\en
\end{lemma}

\begin{proof}
It is clear that if $\rank(H_1(\pi))>1$, then $\pi/\pi_r^{(n)}$ is not cyclic. Now assume that
$\rank(H_1(\pi))=1$. The first statement is immediate. Now let $\phi:\pi\to \Z$ be either of the
two epimorphisms. Clearly $\ker(\phi)=\pi_r^{(1)}$.
 For the second statement, note that by Equation
(\ref{def:deltantwo}) we have
\[ \delta_0(\pi,\phi)=\dim_{\Q}(\pi_r^{(1)}/[\pi_r^{(1)},\pi_r^{(1)}]\otimes_{\Z} \Q).\]
Recall that  $\pi_r^{(1)}/\pi_r^{(2)}$ is torsion free and that
$\pi_r^{(1)}/[\pi_r^{(1)},\pi_r^{(1)}]\otimes_{\Z} \Q=\pi_r^{(1)}/\pi_r^{(2)}\otimes_{\Z} \Q$. This
means that $\delta_0(\pi,\phi)=0$ if and only if $\pi_r^{(1)}=\pi_r^{(2)}$. But in that case
$\pi_r^{(1)}=\pi_r^{(n)}$ for all $n$. The second statement is now clear.
\end{proof}

 Before we state our main theorems we need to introduce admissible triples:

\begin{definition}
Let $\pi$ be a group and $\phi:\pi\to \Z$ be an epimorphism. Let $\varphi_1:\pi \to G_1$ and
$\varphi_2:\pi \to G_2$ be epimorphisms to  locally indicable and amenable groups $G_1$ and $G_2$.
Following \cite[Definition 2.1]{Ha06} we call $(\varphi_1,\varphi_2,\phi)$ an {\em admissible
triple for $\pi$} if there exist epimorphisms $\psi:G_1\to G_2$ and $\phi_2:G_2\to \Z$ such that
$\varphi_2=\psi\circ \varphi_1$, and $\phi=\phi_2\circ \varphi_2$.
\end{definition}

We note that for an admissible triple $(\varphi_1, \varphi_2, \phi)$ for $\pi$ we have the
following commutative diagram and that $(\varphi_i,\phi)$ are admissible pairs for $\pi$ for
$i=1,2$.
 \[   \xymatrix { &G_1\ar[d]^{\psi}\\ \pi \ar[dr]_{\phi} \ar[ur]^{\varphi_1} \ar[r]^{\varphi_2}
&G_2 \ar[d]^{\phi_2}\\& \Z. }
     \]
Clearly for any epimorphism $\phi : \pi \to \Z$ and for any $n \geq m \geq 1$ we  get admissible
triples $(\pi\to \pi/\pi_r^{(n)}, \pi\to \pi/\pi_r^{(m)},\phi)$ for $\pi$.

%===============================================================================\\
\section{The statement of the main results} \label{section:results}

It follows from Theorem \ref{thm:taudelta} and Lemmas \ref{lemma:k1} and \ref{lem:cyclic} that the
following theorem implies Theorem \ref{mainthm1}.

%It follows
%from the universal property of the Cohn localization and from the definition of an admissible
%triple that $\psi:G_1\to G_2$ induces a unique map $C(G_1)\to C(G_2)$ (which we also denote by
%$\psi$) which makes the following diagram commutative:
%\[ \xymatrix{ \Z[G_1] \ar[d]^{\psi} \ar[r]& C(G_1)\ar[d]^{\psi} \\ \Z[G_2] \ar[r] & C(G_2) }
%\]

\begin{theorem} \label{mainthm}
Let $M$ be a  3--manifold with empty or toroidal boundary. Let
$(\varphi_1:\pi_1(M)\to G_1,\varphi_2:\pi_1(M)\to G_2,\phi)$ be an
admissible triple for $\pi_1(M)$. If $\tau(M,\K(G_2')(t))\ne 0$,
then $\tau(M,\K(G_1')(t))\ne 0$ and
\[ \deg(\tau(M,\K(G_1')(t)) = \deg(\tau(M,\K(G_2')(t))+2k \]
for some $k\geq 0$.
\end{theorem}

 The proof of this theorem is postponed to Section \ref{section:mainthm}.\\

Now let $M$ be a 3 manifold. We pick a primitive element $\phi\in H^1(M;\Z)$. Note that if
$b_1(M)=1$ then there is a unique $\phi$ up to sign. Clearly $(\phi,\phi)$ is an admissible pair
and $\K(\Z')(t) = \Q(t)$ and we obtain the corresponding Reidemeister torsion $\tau(M,\Q(t))$.

Given a link $L\subset S^3$ we write $X(L)=S^3\sm N(L)$ for the link exterior, where $N(L)$ denotes
an open tubular neighborhood of $L$. Now we can formulate the following corollary to Theorem
\ref{mainthm}.
\begin{corollary} \label{cor:parity}
 Let $M$ be a 3--manifold and $(\v:\pi_1(M)\to G, \phi)$ an admissible pair for $\pi_1(M)$.
\bn \item If $M=X(K)$ is a knot exterior, then
$\deg(\tau(X(K),\K(G')(t)))$ is odd. \item If $M$ is closed and
$\tau(M,\Q(t))\ne 0$, then $\deg(\tau(M,\K(G')(t)))$ is even.
\item If $M=X(L)$ is the exterior of a link $L=L_1\cup \dots
\cup L_m\subset S^3$ and if $\tau(X(L),\Q(t))\ne 0$, then
\[\deg(\tau(X(L),\K(G')(t)))\equiv \sum_{i=1}^m \phi(\mu_i)\big( 1+\sum_{j \ne i} \lk(L_i,L_j)\big) \modd
2,\] where $\mu_i$ denotes the meridian of $L_i$, $1\le i \le m$.
\en
\end{corollary}

\begin{proof}[{\bf Proof of Corollary \ref{cor:parity}}]
First recall the well--known fact that the Alexander polynomial $\Delta_K(t)$ is of even degree for
any knot $K$, this implies that $\tau(X(K),\Q(t))=\Delta_K(t)(t-1)^{-1}$ is of odd degree. The
first statement now follows immediately from Theorem \ref{mainthm} applied to the admissible triple
$(\v,\phi,\phi)$.

The parity of the multivariable Alexander polynomials for links (cf. \cite[Theorem~1.7.1]{Tu86})
and for closed 3--manifolds (cf. \cite[p.~141]{Tu86}) are well--known. Using standard results (cf.
\cite[Theorem~1.1.2]{Tu86} and also \cite[Theorem~3.4]{FK05}, \cite{FV06}) these can be translated
into parity results for the ordinary one--variable Reidemeister torsion $\tau(M,\Q(t))$. These
results in turn imply the last two statements by Theorem \ref{mainthm}.
\end{proof}

%Before giving the proof of Corollary \ref{cor:parity},
We now discuss a special case of Corollary \ref{cor:parity}(1).
%The proof of Corollary \ref{cor:parity} will be given after the remark below.
Let $K$ be a knot and let $\phi:H_1(X(K);\Z)\to \Z$ be an isomorphism. We write
$\delta_n(K):=\delta_n(\pi_1(X(K)),\phi)$. By Lemma \ref{lemma:k1} this agrees with  Cochran's
\cite{Co04} definition.
 %In \cite{Co04} Cochran
%investigates the degree of \emph{the $n^{\text th}$ order Alexander polynomial} for a knot $K$,
%which is denoted by $\delta_n(K)$. By definition $\delta_n(K)$ is the degree of
%$\Delta_{\phi,1}^\v (X(K))$ where $X(K)$ is the exterior of $K$ in $S^3$, $\phi : \pi_1(X(K))
%\to \Z$ is the abelianization, and $\v : \pi_1(X(K)) \to \pi_1(X(K))/\pi_1(X(K))^{(n+1)}$ is the
%natural projection homomorphism. (For a group $G$, $G^{(n)}$ is inductively defined as follows:
%$G^{(0)} := G$ and $G^{(n+1)} := [G^{(n)},G^{(n)}]$.)
In Question (5) in \cite[Section 14]{Co04} Cochran asks if there
is a knot $K$ and some $n>0$ for which $\delta_n(K)$ is a non-zero
\emph{even} integer. Theorem \ref{thm:taudelta} together with
Corollary \ref{cor:parity}  now gives a negative answer to this
question. More precisely we have the following result.
%It follows from $\delta_n(K) = \deg(\tau(X(K), \K(G')(t)))$ where $G :
%= \pi_1(X(K))/\pi_1(X(K))^{(n+1)}$. Note that $\delta_0(K)$ is the degree of the (classical)
%Alexander polynomial of $K$, hence it is always even. Therefore from Corollary
%\ref{cor:parity}(1) we obtain the following corollary which answers Question (5) in
%\cite[Section 14]{Co04} in the negative.

\begin{theorem}\label{thm:paritydelta}
Let $(n_i)_{i\in \N_0}$ be a sequence of non--negative integers. Then there exists a knot  $K$ with
$\delta_0(K)-1=n_0$ and $\delta_i(K)=n_i$ for $i\geq 1$ if and only if $(n_i)$ is a never
decreasing sequence of odd numbers which is bounded.
\end{theorem}

\begin{proof}
Note that $\delta_0(K)$ is the degree of the (classical) Alexander polynomial of $K$, hence it is
always even. One direction of the theorem now follows from Theorem \ref{thm:taudelta} together with
Corollary \ref{cor:parity} and from the sequence of inequalities (1) in Section
\ref{section:intro}. The realization results in \cite[Theorem~7.3]{Co04} can be strengthened to
prove the converse (cf. also Question (5) in \cite{Co04}).
\end{proof}

\begin{remarks}
\bn \item Note that more realization results are given in
\cite[Section~11]{Ha05}.
\item The statement corresponding  to
Corollary \ref{cor:parity} for the Reidemeister torsion of a 3--manifold associated to a general
linear representation $\pi_1(M) \to \gl(\F,k)$ (where $\F$ is a field) is not known. A careful
analysis of Heegard decompositions of closed 3--manifolds as in \cite{He83} could perhaps be
used to show that for a closed 3--manifold the degrees of twisted Reidemeister torsions are
always even.
 \en
\end{remarks}

Using Corollary \ref{cor:parity} we will prove in Section \ref{section:thurstonparity} the
following theorem which by Lemmas \ref{lemma:k1} and \ref{lem:cyclic} is clearly a slight
generalization of Theorem \ref{mainthm2}.

\begin{theorem}\label{thm:thurstonparity}
Let $M$ be a closed 3--manifold or the exterior of a link in $S^3$. Let $(\v:\pi_1(M)\to G,
\phi)$ be an admissible pair.  If $\tau(M,\Q(t))\ne 0$, then
 \[ \max\{0,\deg(\tau(M,\K(G')(t)))\} \equiv \tnphi \, \modd 2.\]
\end{theorem}

%\noindent The proof of Theorem \ref{thm:thurstonparity} will be given .

%===============================================================================\\
\section{Proof of the theorems}\label{section:proof}

%==================================
\subsection{Cohn localization} \label{section:cohn}

Let $\psi:G_1\to G_2$ be an epimorphism between locally indicable amenable groups $G_1$ and
$G_2$.
 Let $\S$ be the set of all matrices over $\Z[G_1]$ which become invertible under the map
$\Z[G_1]\xrightarrow{\psi} \Z[G_2]\to \K(G_2)$. Denote by $c:\Z[G_1]\to C(\psi)$ the Cohn
localization of $\Z[G_1]$ corresponding to $\S$. Recall that the Cohn localization is characterized
by the following two conditions:
\bn
\item
The homomorphism $c$ is $\S$--inverting, i.e. any matrix in $\S$ becomes invertible under $c$.
\item   The homomorphism $c$  has the universal $\S$--inverting property, i.e. if $d:\Z[G_1]\to D$ is another
$\S$--inverting ring homomorphism, then there exists a unique homomorphism $C(\psi)\to D$ to make
the following diagram commute
\[  \xymatrix {
\Z[G_1] \ar[dr]_{d} \ar[r]^{c} &C(\psi) \ar[d]\\& D. }
     \]
\en
We refer to \cite[Chapter 7]{Co85} for more details. Note that for $r\ne s$ we have $C(\psi)^r
\ncong C(\psi)^s$ since $C(\psi)^i\otimes_{C(\psi)}\K(G_2)\cong \K(G_2)^i$.

\begin{lemma} \label{lem:injinv}
Let $\psi:G_1\to G_2$ be an epimorphism between locally indicable amenable groups. Let $B$ be an
$r\times r$--matrix over $\Z[G_1]$. If $\psi(B):\Z[G_2]^r \to \Z[G_2]^r$ is invertible over
$\K(G_2)$, then $B$ is invertible over $\K(G_1)$.
\end{lemma}

\begin{proof}
First note that  $\psi(B):\Z[G_2]^r\to \Z[G_2]^r$ is injective since $\Z[G_2]\subset \K(G_2)$. Now
let $H=\ker\{\psi:G_1\to G_2\}$. Clearly $H$ is again locally indicable. Note that $B:\Z[G_1]^r\to
\Z[G_1]^r$ can also be viewed as a map between free $\Z[H]$--modules. Pick any right inverse
$\l:G_2\to G_1$ of $\psi$. It is easy to see that  $g\otimes h \mapsto g\l(h)\otimes 1, g\in G_1,
h\in G_2$ induces an isomorphism
\[  \Z[G_1]\otimes_{\Z[G_1]} \Z[G_2]\to \Z[G_1]\otimes_{\Z[H]}\Z. \]
By assumption $B\otimes \id:\Z[G_1]^r\otimes_{\Z[H]} \Z \to \Z[G_1]^r\otimes_{\Z[H]} \Z$ is
injective. Since $H$ is locally indicable it  follows immediately from \cite{Ge83} or
\cite{HS83} (cf. also \cite{St74} for the case of PTFA groups) that $\Z[G_1]^r\to \Z[G_1]^r$ is
injective.

  Since $\K(G_1)$ is flat over $\Z[G_1]$ it follows that
$B:\K(G_1)^r\to \K(G_1)^r$ is injective. But an injective homomorphism between vector spaces
over a skew field of the same dimension is in fact an isomorphism. This shows that $B$ is
invertible over $\K(G_1)$.
\end{proof}

It follows from the universal property of the Cohn localization
and Lemma \ref{lem:injinv} that there exist unique maps
$C(\psi)\to \K(G_1)$ and $ C(\psi)\to \K(G_2)$ which make the
following diagram commute:
 \[   \xymatrix { &\K(G_1)\\ \Z[G_1] \ar[dr]_{\psi} \ar[ur] \ar[r]
&C(\psi) \ar[u]  \ar[d] \\& \K(G_2). }
     \]
Note that since the homomorphism $\Z[G_1]\to \K(G_1)$ is an
injection it follows that the homomorphism $\Z[G_1]\to C(\psi)$ is
an injection as well. We will henceforth identify $\Z[G_1]$ with
its image in $C(\psi)$.

%The following lemma shows that for a 3--manifold $M$ and $\pi :=
%\pi_1(M)$, $\tau(M,C(\psi))\in K_1(C(\psi))/\pm\pi_1(M)$ is
%non--zero.

\begin{lemma} \label{lemmahqg}
Let $M$ be a 3--manifold and $(\varphi_1:\pi_1(M)\to
G_1,\varphi_2:\pi_1(M)\to G_2,\phi)$ an admissible triple. If
$\tau(M,\K(G_2))\ne 0$, then $\tau(M,C(\psi))\ne 0$.
\end{lemma}
\begin{proof}
We have to show that if $H_*(M;\K(G_2))=0$, then $H_*(M;C(\psi))=0$. Denote the chain complex of
the universal cover $\ti{M}$ by $(C_*,d_*)$. Write $\pi:=\pi_1(M)$. By assumption we have
$H_*(M;\K(G_2))=H_*(C_*\otimes_{\Z[\pi]} \K(G_2))=0$. Therefore there exist chain homotopies
$s_k:C_k\otimes_{\Z[\pi]} \K(G_2)\to C_{k+1}\otimes_{\Z[\pi]}  \K(G_2)$ such that for any $k$
\[ d_{k+1} \circ s_k+ s_{k-1}\circ d_k= \id.\]
Note that $C_k \otimes_{\Z[\pi]} \Z[G_2]$ is a free $\Z[G_2]$--module and $C_k\otimes_{\Z[\pi]}
\K(G_2)=(C_k \otimes_{\Z[\pi]} \Z[G_2])\otimes_{\Z[G_2]}\K(G_2)$. Since $\K(G_2)$ is the Ore
localization of $\Z[G_2]$ the chain homotopies $s_k$ induce  $\Z[G_2]$--homomorphisms $s_k':C_k
\otimes_{\Z[\pi]} \Z[G_2] \to C_{k+1}\otimes_{\Z[\pi]} \Z[G_2] $ such that $d_{k+1} \circ s_k'+
s_{k-1}'\circ d_k$ is invertible over $\K(G_2)$.

Furthermore note that $C_k \otimes_{\Z[\pi]} \Z[G_1]$ is a free
$\Z[G_1]$--module and $C_k\otimes_{\Z[\pi]} \Z[G_2]=(C_k
\otimes_{\Z[\pi]} \Z[G_1])\otimes_{\Z[G_1]}\Z[G_2]$, hence we can
find lifts of $s_k'$ to $\Z[G_1]$--homomorphisms $t_k':C_k
\otimes_{\Z[\pi]} \Z[G_1] \to C_{k+1}\otimes_{\Z[\pi]} \Z[G_1] $
such that $d_{k+1} \circ t_k'+ t_{k-1}'\circ d_k$ becomes
invertible under the map $\Z[G_1]\to \Z[G_2] \to \K(G_2)$. From
the definition of the Cohn localization it follows that
\[ d_{k+1} \circ t_k'+ t_{k-1}'\circ d_k :C_k\otimes_{\Z[G_1]}C(\psi)\to C_k\otimes_{\Z[G_1]}C(\psi) \]
is an isomorphism for any $k$, hence $H_*(M;C(\psi))=0$.
\end{proof}

%We can therefore consider the Reidemeister torsion $\tau(M,G)\in K_1(C(G))/\pm \im\{G\to
%C(G)\}$.
%Let $(\v:\pi_1(M)\to G,\phi)$ be an admissible pair for $\pi_1(M)$. Assume that we have
%$\tau(M,\v)\ne 0 \in K_1(C(G))$. Recall that there exists a homomorphism $c: C(G)\to \K(G)$.
%Note that $c$ induces a map $K_1(C(G))\to K_1(\K(G')(t))$. The advantage of $K_1(\K(G')(t)$ is
%that we have a degree function $\deg:K_1(\K(G')(t))\to \Z$.
% We
%define $\deg(\tau(M,\v))$ as the degree $c(\tau(M,\v))\in K_1(\K(G')(t))/\pm \im\{G\to
%\K(G')(t)\}$.

%==================================
\subsection{Proof of Theorem \ref{mainthm}} \label{section:mainthm}
In the following assume that $R$ is a ring equipped with a
(possibly trivial) involution $r\mapsto \ol{r}$ such that
$\ol{r\cdot s}=\ol{s}\cdot \ol{r}$.   This extends to an
involution on $K_1(R)$ via
$\ol{(a_{ij})}:=\left(\ol{a_{ji}}\right)$ for $(a_{ij})\in
K_1(R)$. Note that if $R$ is a skew field, then
$\det(\ol{A})=\ol{\det(A)}$ for any $A\in K_1(R)$.

We need the following lemma.
\begin{lemma} \label{lemma:key}
Let $\k$ be a skew field with (possibly trivial) involution. Let $\ktfield$ be the (skew)
quotient field of a Laurent polynomial ring $\kt$ equipped with the involution given by
$\ol{kt^l}:=t^{-l}\ol{k}$. Assume that $\tau \in K_1(\ktfield)$ has the property that $\tau=kt^j
\, \ol{\tau}$ for some $k\in \K, j\in \Z$. Then
\[  \deg(\tau) \equiv j \mod 2.   \]
\end{lemma}

\begin{proof}
Recall that the Dieudonn\'e determinant defines an isomorphism $K_1(\ktfield)\to
\ktfield^\times/[\ktfield^\times,\ktfield^\times]$. We can therefore rewrite the assumption as
\[ \det(\tau)=kt^j\,  \ol{\det(\tau)} \, \prod_{i=1}^n [f_i(t),g_i(t)] \in \ktfield,\]
for some $f_i(t),g_i(t)\in \ktfield^\times$. For
$p(t)=a_rt^r+a_{r+1}t^{r+1}+\dots+a_st^s \in \kt$ with $a_i\in \K,
a_r\ne 0, a_s\ne 0, r\le s$ we define $l(p(t)):=r$ (the lowest
exponent), $h(p(t)):=s$ (the highest exponent). Clearly $l$ and
$h$ define homomorphisms $\kt\sm \{0\} \to \Z$. Furthermore since
$\ktfield$ is the Ore localization of $\kt$ we can extend $l$ and
$h$ to homomorphisms $l,h:\ktfield^\times\to \Z$ via
$l(f(t)g(t)^{-1})=l(f(t))-l(g(t))$ and
$h(f(t)g(t)^{-1})=h(f(t))-h(g(t))$ for $f(t),g(t)\in \kt\sm
\{0\}$. Note that $l,h$ vanish on commutators. Then we get the
following equalities
\[ \ba{rcl} h(\det(\tau))
&=& h(kt^j\,  \ol{\det(\tau)} \, \prod_{i=1}^n [f_i(t),g_i(t)])\\
&=& j +h(\ol{\det(\tau)})\\
&=& j -l(\det(\tau)).\ea
\]
It follows that
\[ \deg(\tau)\equiv \deg(\det(\tau))\equiv h(\det(\tau))-l(\det(\tau)) \equiv
h(\det(\tau))+l(\det(\tau))\equiv j \mod\, 2.\]
\end{proof}

 We need the following duality
theorem.

\begin{theorem} \label{thm:duality}
Let $M$ be a 3--manifold with empty or toroidal boundary. Assume that $R$ is equipped with a
(possibly trivial) involution  and that $\varphi:\Z[\pi_1(M)]\to R$ is a homomorphism such that
$\v(g^{-1})=\ol{\v(g)}$ for all $g\in \pi_1(M)$.  If $\tau(M,R)\ne 0$, then
\[ \tau(M,R) =\ol{\tau(M,R)} \in K_1(R)/\pm \pi_1(M).\]
\end{theorem}

\begin{proof}
Similarly to $\tau(M,R)$ we can define $\tau(M,\partial M,R)$. By \cite[Theorem 14.1]{Tu01} it
follows that $H_*(M,\partial M;R)=0$ and
\[  \tau(M,R) =\ol{\tau(M,\partial M,R)} \in
K_1(R)/\pm \pi_1(M).\] In particular we are done if $M$ is closed. If $M$ has boundary, then it
follows from the long exact sequence of the pair $(M,\partial M)$ that $H_*(\partial M;R)=0$.
Furthermore by \cite[Theorem~3.4]{Tu01} we have
\[ \tau(M,R)=\tau(M,\partial M,R) \tau(\partial M,R) \in
K_1(R)/\pm \pi_1(M).\] It therefore suffices to show that $\tau(\partial M,R)=1 \in K_1(R)/\pm
\pi_1(M)$. But this follows from an easy argument using the standard CW--structure of a torus.
%Now consider a torus component $T$ of $\partial M$. We equip it
%with a CW--structure with one 0--cell, two 1--cells and one
%2--cell. Denote the generators for $\pi_1(T)$ corresponding to the
%two 1--cells by $a,b$. We get the following chain complex for
%$C_*(\ti{T})$:
%\[ 0\to \Z[a^{\pm 1},b^{\pm 1}]\xrightarrow{\bp a-1\\b-1\ep}\Z[a^{\pm 1},b^{\pm 1}]^2
%\xrightarrow{\bp 1-b & a-1\ep} \Z[a^{\pm 1},b^{\pm 1}]\to 0.\]
%Since $H_*(T;R)=0$ it follows that $\v(a)\ne 0$ or $\v(b)\ne 0$.
%We will assume the latter. It follows that
%$\tau(T,R)=(\v(b-1))(\v(1-b))^{-1})=-1 \in K_1(R)$ (cf.
%\cite[Theorem 2.2]{Tu01} and \cite[Proposition 2.3]{Fr06}).
\end{proof}

For a locally indicable amenable group $G$, we will always consider the group ring $\Z[G]$ together
with the involution given by $\ol{g}:=g^{-1}$ for $g\in G$. Note that this involution extends to
involutions on $\K(G)$ and $K_1(\K(G))$. Recall that for an admissible pair $(\v:\pi \to G, \phi)$
for $\pi$ we have the isomorphism $\Z[G']\tpm \xrightarrow{\cong} \Z[G]$. If we equip $\Z[G']\tpm$
with the involution given by $\ol{gt^l}:=t^{-l}g^{-1}$ for $g\in G'$ then this isomorphism
preserves involutions. The extended involutions on $\K(G')(t)$ and $\K(G)$ are also preserved under
the corresponding isomorphism, and we have a similar result for $K_1(\K(G')(t))$ and $K_1(\K(G))$.

% Given an admissible triple we can pick splittings $\Z\to G_i$ of $\varphi_i,
%i=1,2$ which make the following diagram commute:
% \[   \xymatrix { \Z \ar[r] \ar[dr]& G_1 \ar[d]^{\psi} \\&G_2.  }    \]
%We therefore get an induced commutative diagram of
% ring homomorphisms
% \[   \xymatrix { \Z[\pi_1(M)] \ar[r] \ar[dr]& \Z[G_1']\tpm \ar[d]^{ \psi} \\& \Z[G_2']\tpm. }    \]
% Note that we are suppressing the notation for the
% twisting in the skew Laurent polynomial rings.
%Denote the $\phi$--compatible maps $\Z[\pi_1(M)]\to \K(G_i')\tpm, i=1,2$ by $\varphi_i$ as well.

\begin{proof}[Proof of Theorem \ref{mainthm}]
Let $M$ be a  3--manifold with empty or toroidal boundary. Let $(\varphi_1:\pi_1(M)\to
G_1,\varphi_2:\pi_1(M)\to G_2,\phi)$ be an admissible triple for $\pi_1(M)$. Assume that
$\tau(M,\K(G_2')(t))\ne 0$. Harvey \cite{Ha06} showed that this implies that
$\tau(M,\K(G_1')(t))\ne 0$ and furthermore that $ \deg(\tau(M,\K(G_1')(t)) \geq
\deg(\tau(M,\K(G_2')(t))$. (See also \cite[Theorem 1.3]{Fr06}, we also refer to \cite{Co04} in the
case of a knot exterior). Therefore in order  to prove Theorem \ref{mainthm} it is enough to show
that
\[  \deg(\tau(M,\K(G_1')(t))  \equiv \deg(\tau(M,\K(G_2')(t)) \mod 2.\]

 Write $\pi:=\pi_1(M)$. Note that $\tau(M,C(\psi))$ and $\ol{\tau(M,C(\psi))}$ are
non--zero in $K_1(C(\psi))/\pm \pi \cup \{0\}$ by Lemma \ref{lemmahqg}. Now pick a representative
$\tau^{rep}(M,C(\psi))$ of $\tau(M,C(\psi))$ in $K_1(C(\psi))$. Note that
$\ol{\tau^{rep}(M,C(\psi))}$ is a representative of $\ol{\tau(M,C(\psi))}$ in $K_1(C(\psi))$. It
follows from Theorem \ref{thm:duality} that
\[ \tau^{rep}(M,C(\psi))=\eps g_1\ol{\tau^{rep}(M,C(\psi))} \in K_1(C(\psi)) \]
for some $\eps \in \{1,-1\}, g_1\in G_1$. (Recall that we identified $\Z[G_1]$ with the image of
$\Z[G_1]$ in $C(\psi)$.) Recall that we have homomorphisms $\a_i:C(\psi)\to \K(G_i')(t)$ for
$i=1,2$ induced from the universal property of the Cohn localization $C(\psi)$. These homomorphisms
induce homomorphisms $\a_i:K_1(C(\psi))\to K_1(\K(G_i')(t))$ for $i=1,2$.
%The following diagram can be used for orientation:
%\[ \xymatrix { &\Z[G_1]\ar[r]\ar[d]^\psi & C(G_1)\ar[r]^{c_1}\ar[d]^\psi &\K(G_1)=\K(G_1')(t)\\
%\Z[\pi] \ar[ur]^{\v_1} \ar[r]^{\v_2} \ar[dr]^{\phi}&\Z[G_2] \ar[d] \ar[r]&C(G_2)\ar[d]\ar[r]^{c_2}&\K(G_2)=\K(G_2')(t)\\
%&\zt\ar[r]&\Q(t). } \]
Then for $i=1,2$,
\[\a_i(\tau^{rep}(M,C(\psi)))=\a_i(\eps g_1\ol{\tau^{rep}(M,C(\psi)}))=
\eps \a_i(g_1)\a_i(\ol{\tau^{rep}(M,C(\psi)}))\in
K_1(\K(G_i')(t)).\]

Note that because Reidemeister torsion is functorial (cf. e.g. \cite[Proposition~3.6]{Tu01}),
$\a_i(\tau^{rep}(M,C(\psi)))\in K_1(\K(G_i')(t))$ is a representative of $\tau(M,\K(G_i')(t))\in
K_1(\K(G_i')(t))/\pm\pi$. Now write $\a_i(g_1)=k_i t^{l_i}$ for some $k_i\in \K(G_i')$ and
$l_i\in \Z$. Note that  $l_i=\phi(g_1)$ since $\Z[\pi]\to \K(G_i')\tpm$ is $\phi$--compatible.
In particular $l_1=l_2$. The theorem now follows from Lemma \ref{lemma:key} which states that
\[   l_i \equiv \deg(\tau(M,\K(G_i')(t))) \mod 2\]
for $i=1,2$.
%Now note that  $c(g_1)=k_1t^l$ for some $k_1\in \K(G_1')$ since
%$\Z[\pi]\to \K(G_1')\tpm$ is $\phi$--compatible. Applying the same
%argument as above for $c$ instead of $\psi$ we get
% \[  \deg(\tau(M,\K(G_1')(t))) \equiv l \mod 2.\]
\end{proof}

%==================================
\subsection{Proof of Theorem
\ref{thm:thurstonparity}}\label{section:thurstonparity}

We will need the following lemma.

\begin{lemma} \label{lemmaeven}
Let $\phi \in H^1(M)$ be primitive such that $\tau(M,\Q(t))\ne 0$. If $M$ is closed, then
$||\phi||_T$ is even. Assume that $\partial M$ consists of a non--empty collection of tori $N_1\cup
\dots \cup N_m$. If $\phi|_{H_1(N_i)}=0$ then let $n_i:=0$, otherwise define
\[ n_i=\max\{ n\in \N | \phi|_{H_1(N_i)}=n h \mbox{ for some }h\in \hom(H_1(N_i),\Z)\},\]
i.e. $n_i$ is the divisibility of $\phi|_{H_1(N_i)}$. Then either $||\phi||_T=0$ or
\[ ||\phi||_T\equiv \bigg(\sum_{i=1}^m n_i \bigg) \modd \, 2.\]
\end{lemma}

\begin{proof}
If $M$ is closed then $S$ is closed, hence $\chi(S)$ is even. Now assume that $\partial M$ is a
non--empty collection of tori. Since $\tau(M,\Q(t))\ne 0$ it follows from \cite{Mc02} that there
exists a connected Thurston norm minimizing surface $S$ dual to $\phi$. Assume that $||\phi||_T>0$,
then $||\phi||_T=-\chi(S)$. Note that
\[ -\chi(S)\equiv b_0(\partial S) \,\,  \modd 2,\]
this follows from the observation that adding a 2--disk to each component of $\partial S$ gives
a closed surface, which has even Euler characteristic. Now consider $N_i$. Clearly $S\cap N_i$
is Poincar\'e dual to $\phi|_{H_1(N_i)}$. It follows from a standard argument that, modulo 2,
$\partial S \cap N_i$ has $n_i$ components.
%\footnote{I slightly had to change the lemma, the
%old version proved that $\chi(S)\equiv  \bigg(\sum_{i=1}^s n_i \bigg) \modd \, 2$, which is not
%quite the same, since $\chi(S)\ne ||\phi||_T$ in general. For example the old lemma was
%incorrect for $S^1\times D^2$.}
\end{proof}

Now we can give the proof of Theorem \ref{thm:thurstonparity}.

\begin{proof}[Proof of Theorem \ref{thm:thurstonparity}]
First assume that $||\phi||_T=0$. Note that it follows from Equation (\ref{equ:lowerbound}) and
from Theorem \ref{thm:taudelta} that  $\deg(\tau(M,\K(G')(t)))\leq ||\phi||_T$.  The required
equality is now immediate. Now assume that $||\phi||_T>0$.

Since both sides are $\N$--linear we can restrict ourselves to the case that $\phi$ is
primitive. If $M$ is closed, then $||\phi||_T$ is even and the theorem follows from Corollary
\ref{cor:parity}(2). Now assume that $M=X(L)$ where $L=L_1\cup \dots \cup L_m$ is a link in
$S^3$.

Denote by $\mu_1,\dots,\mu_m$ (respectively $\l_1,\dots,\l_m$) the meridians (respectively
longitudes) of $L_1\cup \dots \cup L_m$. Denote by $N_1\cup \dots \cup N_m$ the boundary components
of $X(L)$. Let $n_1,\dots,n_m$ as in Lemma \ref{lemmaeven}.  Then by Lemma \ref{lemmaeven}
\[ ||\phi||_T\equiv \bigg(\sum_{i=1}^m n_i \bigg) \modd \, 2.\]
Clearly $\mu_i,\l_i$ is a basis for $H_1(N_i)$ and $\phi|_{H_1(N_i)}\in H^1(N_i)$ is Poincar\'e
dual to $\phi(\mu_i)\mu_i+\big(\sum_{j\ne i}\phi(\mu_j)\lk(L_i,L_j)\big)\l_i$. In particular
\[ n_i =\gcd\Big(\phi(\mu_i),\sum_{j\ne i}\phi(\mu_j)\lk(L_i,L_j)\Big)\]
where we set $\gcd(a,0):=\gcd(0,a):=0$.

Without loss of generality we can assume that there exists an $r$ such that $\phi(\mu_i)\equiv  1
\modd 2$ for $i\leq r$ and $\phi(\mu_i)\equiv  0 \modd 2$ for $i> r$. It follows immediately that
\[ \sum_{i=1}^m n_i \equiv r+\sum_{i=r+1}^m \sum_{j\ne
i}\phi(\mu_j)\lk(L_i,L_j))\, \, \modd 2.\]
 On the other hand by Corollary \ref{cor:parity} and by the symmetry of the linking numbers we
 have
\[ \ba{rcll}\deg(\tau(M,\K(G')(t)))&\equiv &\sum_{i=1}^m \phi(\mu_i)+\sum_{i=1}^m \sum_{j \ne i}\phi(\mu_i) \lk(L_i,L_j)& \modd
2\\[2mm]
&\equiv &r+\sum_{i=1}^m \sum_{j< i} (\phi(\mu_i)+\phi(\mu_j))\lk(L_i,L_j)& \modd 2\\[2mm]
&\equiv &r+\sum_{i=r+1}^m \sum_{j< i}(\phi(\mu_i)+\phi(\mu_j)) \lk(L_i,L_j)& \modd 2\\[2mm]
&\equiv &r+\sum_{i=r+1}^m \sum_{j< i} \phi(\mu_j)\lk(L_i,L_j)& \modd 2\\[2mm]
&\equiv &r+\sum_{i=r+1}^m \sum_{j\ne i} \phi(\mu_j)\lk(L_i,L_j)& \modd 2.\ea\]
 Note that we used that $\phi(\mu_i)+\phi(\mu_j)\equiv 0 \modd 2$ for $i,j\leq r$ and
$\phi(\mu_j)\equiv 0 \modd 2$ for $j>r$.
\end{proof}

%\begin{remark}
%Our results here and the monotonicity results in \cite{Co04,Ha06,Fr06} are proved using a
%`bottom--up' approach. More precisely, given an admissible triple and given information on
%$\tau(M,\K(G_2')(t))$ we can prove a statement about $\tau(M,\K(G_1')(t))$.
%
%It is helpful to think about this subject from a 'top--down' perspective. More precisely, our
%philosophy is that given $(M,\phi)$ (at least in the knot complement case) there should exist a
%Laurent polynomial ring $R\tpm$ and an invariant $\Delta_\phi(t)\in R\tpm$
%% (e.g. in $K_1(C(\pi_1(M))$, for some Cohn localization of $C(\pi_1(M))$, or perhaps $\U(G)$ instead)
%which has the following properties:
% \bn
% \item  $\Delta_\phi(t)=0$ or the degree of $\Delta_\phi(t)$
%%(which has to be defined using some so far unknown algebraic tool)
%gives the Thurston norm.
% \item All other higher--order Alexander polynomials or
%Reidemeister torsions are just specializations of $
%\Delta_\phi(t)$ under ring homomorphisms $R\tpm \to \K(G')\tpm$.
%\item $\Delta_\phi(t)$ is symmetric in an appropriate sense. \en
% Property (2)  would explain why
%the degrees can only go down, furthermore property (3) would
%explain why the degrees go down only by even numbers. Together
%with (1) this would also explain Theorem \ref{thm:thurstonparity}.
%\end{remark}

\end{document}